\documentclass{amsart}
\usepackage{amssymb}
\usepackage{verbatim}
\begin{document}
\newcommand\im{\operatorname{Im}}
\newcommand\re{\operatorname{Re}}
\newcommand\Id{\operatorname{Id}}
\newcommand\Real{\mathbb{R}}
\newcommand\Cinf{{\mathcal C}^{\infty}}
\newcommand\dist{{\mathcal C}^{-\infty}}
\newcommand\Diff{\operatorname{Diff}}
\newcommand\pa{\partial}
\newcommand\supp{\operatorname{supp}}

\setcounter{secnumdepth}{3}
\newtheorem{lemma}{Lemma}[section]
\newtheorem{prop}[lemma]{Proposition}
\newtheorem{thm}[lemma]{Theorem}
\newtheorem{cor}[lemma]{Corollary}
\newtheorem{result}[lemma]{Result}
\newtheorem*{thm*}{Theorem}
\newtheorem*{prop*}{Proposition}
\newtheorem*{cor*}{Corollary}
\newtheorem*{conj*}{Conjecture}
\numberwithin{equation}{section}
\theoremstyle{remark}
\newtheorem{rem}[lemma]{Remark}
\newtheorem*{rem*}{Remark}
\theoremstyle{definition}
\newtheorem{Def}[lemma]{Definition}
\newtheorem*{Def*}{Definition}
\renewcommand{\theenumi}{\roman{enumi}}
\renewcommand{\labelenumi}{(\theenumi)}

\title[Absence of super-exponentially decaying eigenfunctions]
{Absence of super-exponentially decaying eigenfunctions
on Riemannian manifolds with pinched negative curvature}
\author[Andras Vasy]{Andr\'as Vasy}
\author{Jared Wunsch}
\address{Department of Mathematics, MIT and Northwestern University}
\address{Department of Mathematics, Northwestern University}
\email{andras@math.mit.edu}
\email{jwunsch@math.northwestern.edu}
\date{October 21, 2004}
\thanks{A.\,V. was partially supported by NSF grant DMS-0201092, a Clay 
Research Fellowship and a Fellowship from the Alfred P.\,Sloan Foundation.}
\thanks{J.\,W. was partially supported by NSF grants DMS-0323021 and
DMS-0401323.}
\maketitle

\section{Introduction and statement of results}

Let $(X,g)$ be a metrically complete, simply connected Riemannian
manifold with bounded geometry and pinched negative curvature,
i.e.\ there are constants $a>b>0$ such that $-a^2<K<-b^2$ for all
sectional curvatures $K$. Here bounded geometry is used in the
sense of Shubin, \cite[Appendix~1]{Shubin:Spectral}, namely
that all covariant derivatives of the Riemannian curvature tensor are
bounded and the injectivity radius is uniformly bounded below by a positive
constant. We show that there are no superexponentially
decaying eigenfunctions of $\Delta$ on $X$; here $\Delta$ is the positive
Laplacian of $g$.
That is, fix some $o\in X$, and let $r(p)=d(p,o)$,
$p\in X$, be the distance function. Then:

\begin{thm}\label{thm:Lap}
Suppose that $(X,g)$ is as above. If $(\Delta-\lambda)\psi=0$
and
$\psi\in e^{-\alpha r}L^2(X)$ for all $\alpha$, then $\psi$ is identically $0$.
\end{thm}

Since the curvature assumptions imply exponential volume growth, and due to
elliptic regularity, the $L^2$ norm may be replaced by any $L^p$, or
indeed Sobolev, norm. This result strengthens Mazzeo's unique
continuation theorem at infinity \cite{Mazzeo:Unique} by eliminating
the asymptotic curvature assumption (4) there.

As shown below, the negative curvature assumption enters via the strict
uniform convexity of the geodesic spheres centered at $o$, much as
in the work of Mazzeo \cite{Mazzeo:Unique}. Thus, as observed by Rafe Mazzeo,
the arguments go through equally well if $X$ is replaced by a manifold $M$
which is the union of a `core' $M_0$ (not necessarily compact) and a product
manifold $M_1=(1,\infty)_r\times N$, with a Riemannian metric
$g=dr^2+k(r,.)$, $k$ a metric on $N(r)=\{r\}\times N$, $M_1$ having bounded
geometry,
provided that the second fundamental
form of $N(r)$ is strictly positive, uniformly in $r$. Indeed, the assumptions
on $\psi$ only need to be imposed on $M_1$, see Remark~\ref{rem:loc}.

Following \cite[Appendix~1]{Shubin:Spectral}
we remark that an equivalent formulation of the definition
of a manifold of bounded geometry is the requirement that the injectivity
radius is bounded below by a positive constant $r_{\mathrm{inj}}$,
and that the transition functions between intersecting
geodesic normal coordinate charts (called canonical coordinates in
\cite{Shubin:Spectral})
of radius $<r_{\mathrm{inj}}/2$, say, are $\Cinf$ with uniformly bounded
derivatives (with bound independent of the base points).
We let $\Diff(X)$ denote the algebra of differential operators
corresponding
to the bounded geometry, called the algebra of $\Cinf$-bounded 
differential operators in \cite[Appendix~1]{Shubin:Spectral}. Thus,
in any canonical coordinates, $A\in\Diff^m(X)$ has the form
$\sum_{|\alpha|\leq m} a_\alpha(x)D_x^\alpha$, with $\pa^\beta a_\alpha$
uniformly bounded, with bound $C_{|\beta|}$
independent of the canonical coordinate
chart, for all multiindices $\beta$. We also write $H^k(X)$ for the $L^2$-based
Sobolev spaces below.

If $E$ is a vector bundle of bounded geometry,
in the sense of \cite[Appendix~1]{Shubin:Spectral}, then
Theorem~\ref{thm:Lap} is also valid for $\Delta$ replaced by
any second order
differential operator $P\in \Diff^2(X,E)$
acting on sections of $E$ with scalar principal symbol
equal to that of $\Delta$, i.e.\ the metric function on $T^*X$. Indeed,
we can even localize at infinity, i.e.\ assume $Pu=0$ only near infinity,
and obtain the conclusion that $u$ is $0$ near infinity.

\begin{thm}\label{thm:perturb}
Let $(X,g)$ be as above.
Suppose $P\in\Diff^2(X,E)$, $\sigma_2(P)=g\,\Id$, where $g$ denotes
the metric function on $T^*X$. If $P\psi\in\Cinf_c(X,E)$
and $\psi\in e^{-\alpha r}L^2(X,E)$
for all $\alpha$ then $\psi\in\Cinf_c(X,E)$.
\end{thm}

\begin{rem}
This theorem, together with the standard unique continuation result,
\cite[Theorem~17.2.1]{Hor}, implies that if $P\psi=0$ on $X$ then $\psi=0$
on $X$, just as in Theorem~\ref{thm:Lap}.

In addition, the Sobolev order of the assumptions on $\psi$ and $P\psi$ is
immaterial. In fact, as discussed in Remark~\ref{rem:loc}, the argument
localizes near infinity, hence we may assume $P\psi\in\dist_c(X,E)$,
and then elliptic regularity allows us to conclude that if
$\psi$ is in an exponentially weighted Sobolev space near infinity
(with possibly a negative exponent) then it is in the corresponding
weighted $L^2$-space.
\end{rem}

We also remark on the bounded geometry hypotheses, more precisely on the
assumptions on covariant derivatives. The results of Anderson and Schoen
\cite{Anderson-Schoen:Positive} on harmonic functions on negatively
curved spaces are results below the continuous spectrum. Thus, these are
elliptic problems even at infinity, in a rather strong sense -- stronger
than just the uniform ellipticity on manifolds with bounded geometry
discussed below. In
particular, the notion of positivity and the maximum principle are available,
and can be used to eliminate conditions on covariant derivatives.
However, for eigenfunctions embedded in the continuous spectrum such tools
are unavailable. Indeed, for $\lambda$ large,
$\Delta-\lambda$ can be seen to lose `strong'
ellipticity (so e.g.\ it is not Fredholm on $L^2(X)$),
and is in many ways (micro)hyperbolic, at least in settings
with an additional structure (cf.\ the
discussion in \cite{RBMSpec} in asymptotically flat spaces).
In such a setting commutator estimates
are very natural, and have a long tradition in PDEs; this explains
the role of the assumption on the covariant derivatives.

We are very grateful to Rafe Mazzeo and Richard Melrose for numerous very
helpful conversations, and for their interest in the present work.

\section{The proofs}
To make the argument more transparent, we write up the
proof of Theorem~\ref{thm:Lap}, at each step pointing out any significant
changes that are needed to prove Theorem~\ref{thm:perturb}. The proofs
are a version of Carleman estimates, see \eqref{eq:Carleman} below, at
least for self-adjoint operators, but we phrase these somewhat
differently, in the spirit of operators with complex symbol and codimension
2 characteristic variety (which in this case is in the semiclassical limit)
on which the Poisson bracket of the
real and imaginary part of the (in this case, semiclassical)
principal symbol is positive. This corresponds to non-solvability
of the inhomogeneous PDE in the sense of \cite[Section~26.4]{Hor}; see also
\cite{Zworski:Numerical} for a recent discussion.

We consider eigenfunctions of $\Delta$ that are superexponentially decaying:
$(\Delta-\lambda)\psi=0$, and $\psi\in e^{-\alpha r} L^2(X)$ for all
$\alpha$, $\|\psi\|_{L^2(X)}=1$
(for convenience). Note that $\lambda$ is real by the self-adjointness of
$\Delta$. Moreover, $\psi\in\Cinf(X)$
by standard elliptic regularity, and indeed $\psi\in e^{-\alpha r'}H^m(X)$
for all $m$, where $r'$ is a smoothed version of $r$, changed only near
the origin. It is convenient to assume that $r'>0$, so $\inf_X r'>0$.

For $\alpha$ real, we consider
\begin{equation*}
P_\alpha=e^{\alpha r'}(\Delta-\lambda)e^{-\alpha r'}.
\end{equation*}
Here we need to use $r'$ since $r$ is not smooth at $o$. However, for
notational simplicity, to avoid an additional compactly supported error
term on almost every line, we ignore this, and simply add back a compactly
supported error term in \eqref{eq:h-comm-15}.

Let
\begin{equation*}
\re P_\alpha=\frac{1}{2}(P_\alpha+P_\alpha^*),
\ \im P_\alpha=\frac{1}{2i}(P_\alpha-P_\alpha^*),
\end{equation*}
be the symmetric and skew-symmetric parts of $P_\alpha$. Thus,
$P_\alpha
=\re P_\alpha+i\im P_\alpha$, and $\re P_\alpha$, $\im P_\alpha$ are
symmetric. Note also that $P_\alpha \psi_\alpha=0$ where $\psi_\alpha
=e^{\alpha r}\psi$.
Thus,
\begin{equation}\label{eq:h-comm-8}
0=\|P_\alpha \psi_\alpha\|^2=\|\re P_\alpha \psi_\alpha\|^2+\|\im P_\alpha
\psi_\alpha\|^2
+\langle i[\re P_\alpha,\im P_\alpha]\psi_\alpha,\psi_\alpha\rangle.
\end{equation}
Roughly speaking,
this will give a contradiction provided the commutator is positive
-- although in the presence of error terms one needs to be a little
more careful.

We remark that this argument parallels the last part
of the $N$-body argument of \cite{Vasy:Exponential}, showing exponential
decay and unique continuation results
for $N$-particle Hamiltonians with second order interactions, which in
turn placed the work of Froese and Herbst \cite{FroExp} in potential
scattering into this framework. However, in \cite{Vasy:Exponential}
(as in \cite{FroExp}) this
is the simplest part of the argument; it is much more work to show that
$L^2$-eigenfunctions decay at a rate given by the next threshold
above the eigenvalue $\lambda$ -- hence superexponentially in the absence
of such thresholds.

We now relate our arguments to the usual
Carleman-type arguments, at least if $P_0$ is symmetric (as is for
self-adjoint operators in $\Diff^2(X)$ with the same principal
symbol as $\Delta$).
In those, one considers $P_\alpha$ and $P_{-\alpha}$,
with the same notation as above, and computes $\|P_\alpha \psi_\alpha\|^2\pm
\|P_{-\alpha}\psi_\alpha\|^2$. Since $P_0$ is symmetric, indeed self-adjoint,
$P_{-\alpha}=P_\alpha^*$, so
\begin{equation}\begin{split}\label{eq:Carleman}
&\|P_\alpha \psi_\alpha\|^2+
\|P_{-\alpha}\psi_\alpha\|^2=2\|\re P_\alpha\psi_\alpha\|^2
+2\|\im P_\alpha\psi_\alpha\|^2,\\
&\|P_\alpha \psi_\alpha\|^2-
\|P_{-\alpha}\psi_\alpha\|^2=2\langle i[\re P_\alpha,\im P_\alpha]\psi_\alpha,
\psi_\alpha\rangle.
\end{split}\end{equation}
Thus, the usual Carleman argument breaks up
\eqref{eq:h-comm-8} into two pieces, and is completely equivalent to
\eqref{eq:h-comm-8}. However, dividing up $P_\alpha$ into its symmetric
and skew-symmetric parts makes the calculations below more systematic,
which is particularly apparent in how the double commutator appears
in $\re P_\alpha$ below. This double commutator, in turn, makes it
clear why various terms, which one might expect by expanding out
the squares $\|P_{\pm\alpha} \psi_\alpha\|^2$, do not appear in
the evaluation of $\|P_\alpha \psi_\alpha\|^2\pm
\|P_{-\alpha}\psi_\alpha\|^2$.

Due to the prominent role played by $r$, we work in
Riemannian normal coordinates. So
let $g=dr^2+k(r,.)$ be the metric on $X$, where $k$
is the metric on the geodesic
sphere of radius $r$, denoted by $S(r)$, and let $A(r,.)\,dr\wedge \omega$
denote the volume element, $\omega$ being the standard volume form on the
unit sphere. By the bounded geometry assumptions, $\pa_r\log A
=-\Delta r\in\Cinf_{\mathrm{b}}(X)=\Diff^0(X)$
(see e.g.\ \cite[Lemma~2.3]{Zhu:Comparison}
for the identity),
i.e.\ is uniformly bounded with analogous
conditions on the covariant derivatives.
Then
\begin{equation}\label{eq:Lap-normal}
-\Delta=\pa_r^2+(\pa_r A)\pa_r-\Delta_{S(r)}.
\end{equation}

Now,
\begin{equation*}\begin{split}
&P_\alpha=\Delta-\lambda+e^{\alpha r}[\Delta,e^{-\alpha r}],\\
&\re P_\alpha=\Delta-\lambda+\frac{1}{2}[e^{\alpha r},[\Delta,e^{-\alpha r}]]\\
&\im P_\alpha=\frac{1}{2i}(e^{\alpha r}[\Delta,e^{-\alpha r}]
+[\Delta,e^{-\alpha r}]e^{\alpha r}).
\end{split}\end{equation*}
Here the expressions for $\re P_\alpha$ and $\im P_\alpha$ follow
directly from the definition of the symmetric and skew-symmetric parts,
using that $\Delta$ and $e^{\pm\alpha r}$ are symmetric.

In the double commutator in the expression for $\re P_\alpha$ above, changing
$\Delta$ by a first order operator would not alter the result, as commutation
with a scalar reduces the order by $1$. Thus,
in view of \eqref{eq:Lap-normal},
in the double commutator in $\re P_\alpha$ all terms but $\pa_r^2$ give
vanishing contribution, so we immediately see that
\begin{equation*}
\re P_\alpha=\Delta-\lambda-\alpha^2.
\end{equation*}
We next compute the skew-symmetric part. This is
\begin{equation*}
\im P_\alpha=\frac{1}{i}(2\alpha\pa_r+\alpha(\pa_r\log A)).
\end{equation*}
Thus,
\begin{equation*}
i[\re P_\alpha,\im P_\alpha]=\alpha[\Delta,2\pa_r+(\pa_r\log A)].
\end{equation*}

The crucial estimate for this commutator that we need below is
that there is $c>0$ such that
\begin{equation}\label{eq:comm-8}
[\Delta,2\pa_r+(\pa_r\log A)]\geq c\Delta_{S(r)}+R,\ R\in\Diff^1(X);
\end{equation}
here $R$ is symmetric and the inequality is understood in the sense
of quadratic forms, e.g.\ with domain $H^2(X)$.
Since the commutator is a priori in $\Diff^2(X)$, this means that we merely
need to calculate its principal symbol, which in turn only
depends on the principal symbols of the commutants. Thus,
with $H_g$ denoting the Hamilton
vector field of $g$, and $\sigma$ the canonical dual variable of $r$,
with respect to the product decomposition $(0,\infty)\times S$ of
$X\setminus o$,
the principal symbol of $\pa_r$ is $\sigma_1(\pa_r)=i\sigma$, and
\begin{equation*}
\sigma_2([\Delta,2\pa_r])=2H_g \sigma.
\end{equation*}
It is convenient to
rephrase this by noting that $2\pa_r=-[\Delta,r]+R'$,
$R'\in\Diff^0(X)$, so $2i\sigma=\sigma_1(2\pa_r)=iH_g r$, and hence
$\sigma_2([\Delta,2\pa_r])=H_g^2 r$. The estimate we need then is that
there is $c>0$ such that
\begin{equation}\label{eq:Hessian-est}
H_g^2 r\geq ck.
\end{equation}
Indeed, \eqref{eq:Hessian-est} implies \eqref{eq:comm-8}, since
for each $x\in X$, both sides of \eqref{eq:Hessian-est}
are quadratic forms on $T^*X$, depending smoothly on $x$, so their difference
can be written as $\sum a_{ij}(x)\xi_i\xi_j$ ($\xi_i$ are canonical
dual variables of local coordinates $x_i$), with $a_{ij}$ a non-negative
matrix. This in turn is the principal
symbol of $\sum_{ij}D_{x_i}^*a_{ij}(x)D_{x_j}$, and
\begin{equation*}
\langle \sum_{ij}D_{x_i}^*a_{ij}(x)D_{x_j} v,v\rangle=
\int_X \sum_{ij}a_{ij}(x) D_{x_i}v\,\overline{D_{x_j}v}\,dg\geq 0.
\end{equation*}

To analyze \eqref{eq:Hessian-est}, recall
that arclength parameterized geodesics of $g$ are projections
to $X$ of the integral curves of $\frac{1}{2}H_g$ inside $S^*X$, the unit
cosphere bundle of $X$.
Thus, \eqref{eq:Hessian-est} tells us that $r$ is strictly convex along
geodesics tangent to $S(r_0)$ at the point of contact. Equivalently,
the Hessian $\nabla dr$, which is the form on the fibers of $T X$
dual to $H_g^2 r$, is strictly positive on
$TS(r_0)$, uniformly as $r_0\to\infty$. As $r=r_0$ defines $S(r_0)$, with
$|\nabla r|=1$, this Hessian equals the second fundamental form of $S(r_0)$;
hence \eqref{eq:Hessian-est} is also equivalent to the uniform convexity of
the hypersurfaces $S(r_0)$.

Now, \eqref{eq:Hessian-est} follows immediately when the sectional
curvatures of $X$
are bounded above by a negative constant $-b^2$, since by
the Hessian comparison theorem
(see e.g.\ \cite[Theorem~1.1]{Schoen-Yau:DG}),
$H_g^2 r|_{T_{S(r_0)}X}\geq H_{g_0}^2 r|_{T_{S(r_0)}X}$,
where $g_0$ is the metric with constant negative sectional curvature $-b^2$,
and the right hand side is $b\coth br\geq b$
(cf.\ \cite[Equation~(1.7)]{Schoen-Yau:DG}).

We will consider
$\alpha\to\infty$, but for notational reasons it is convenient to work
in the semiclassical setting.
Thus, let $h=\alpha^{-1}$, $h\in(0,1]$,
$\Delta_h=h^2\Delta$, $\Delta_{S(r),h}=h^2\Delta_{S(r)}$,
and slightly abuse notation by writing
$P_h=h^2e^{r/h}(\Delta-\lambda)e^{-r/h}$, so
\begin{equation*}\begin{split}
&\re P_h=\Delta_h-1-h^2\lambda,\ \im P_h=\frac{1}{i}(2h\pa_r+h(\pa_r\log A)),\\
&i[\re P_h,\im P_h]\geq ch\Delta_{S(r),h}+h^3 R,\ R\in\Diff^1(X).
\end{split}\end{equation*}
We denote the space of semiclassical
differential operators of order $m$ by $\Diff^m_h(X)$. We recall
that $A\in\Diff^m_h(X)$ means that, in the usual multiindex notation,
$A=\sum_{|\alpha|\leq m}
a_\alpha(x) (hD_x)^\alpha$ locally; in our bounded geometry setting
we still impose, as for standard differential operators, that
for all multiindices $\beta$,
$\pa^\beta a_\alpha$ is bounded uniformly in all
Riemannian normal coordinate charts of radius $R$ ($R$ less than half
the injectivity radius, say), with bound only dependent on $|\beta|$.
Then, weakening the above statements somewhat, in a way that still suffices
below,
\begin{equation}\label{eq:h-comm-5}
\re P_h=\Delta_h-1+hR_1',\ \im P_h=\frac{1}{i}(2h\pa_r+hR_2'),
\ i[\re P_h,\im P_h]\geq ch\Delta_{S(r),h}+h^2 R_3'.
\end{equation}
with $R_1',R_2'\in\Diff^1_h(X)$, $R_3'\in\Diff^2_h(X)$.

We stated \eqref{eq:h-comm-5} in a weakened form to make it only depend
on the principal symbol of $\Delta$. Namely,
if $\Delta$ is replaced by any operator $\Delta+Q$, $Q\in\Diff^1(X,E)$
(not necessarily symmetric), and $P'_h=h^2 e^{r/h}(\Delta+Q-\lambda)e^{-r/h}$,
then $P'_h-P_h=h^2 e^{r/h}Qe^{-r/h}\in h\Diff^1_h(X,E)$, so
$\re P'_h-\re P_h,\im P'_h-\im P_h\in h\Diff^1_h(X,E)$, and thus
\begin{equation*}\begin{split}
&i[\re P'_h,\im P'_h]-i[\re P_h,\im P_h]\\
&=i[\re P'_h-\re P_h,\im P'_h]+i[\re P_h,\im P'_h-\im P_h]\in h^2\Diff^2_h(X,E),
\end{split}\end{equation*}
where we used that $\re P_h,\im P_h$ have scalar principal symbols,
hence so do $\re P'_h$ and $\im P'_h$, giving the extra $h$ (compared
to the order of the product) and the lower order in the commutators.
In other
words, \eqref{eq:h-comm-5} still holds for $P_h$ replaced by $P'_h$.

In the above calculations we ignored a compact subset of $X$, so we
need to add a compactly supported error. To avoid overburdening the
notation, we write $r$ for the smoothed out distance function, denoted by
$r'$ above, so for $r$ sufficiently large, $r(p)=d(p,o)$.
Thus, we have shown that
for some $c>0$,
\begin{equation}\label{eq:h-comm-15}
i[\re P_h,\im P_h]\geq ch\Delta_{S(r),h}+h^2R_0+h R_0',\ R_0,R_0'
\in\Diff^2_h(X),
\end{equation}
$R_0'$ supported in $r\leq r_1$ for some $r_1>0$, with the inequality
holding in the sense of
operators. Since $\Delta_{S,h}=\Delta_h-(h D_r)^2-h(D_r \log A) hD_r$,
this estimate implies
\begin{equation}\label{eq:h-comm-16}
\langle i[\re P_h,\im P_h]\psi_h,\psi_h\rangle
\geq\langle (ch+h R_1 \re P_h+h R_2\im P_h+h^2R_3+hR_4)
\psi_h,\psi_h\rangle
\end{equation}
with $R_1\in\Diff^0_h(X)$, $R_2\in\Diff^1_h(X)$ and
$R_3,R_4\in\Diff^2_h(X)$, $R_4$ having compact support in $r\leq r_1$.
(In fact, for our purposes the compact support assumption
is equivalent to assuming that $R_4$ is $o(1)$ as
$r\to\infty$, as it can be absorbed in the first term for $r$
sufficiently large.)

We now show how to use \eqref{eq:h-comm-16} to prove unique
continuation at infinity. To be systematic, we set this part up
somewhat abstractly. Recall that $P_h\in\Diff^2_h(X,E)$ is elliptic
(or more precisely uniformly elliptic, both in $X$ and in $h$)
if there is $C>0$ such that for all $(x,\xi)\in T^*X\setminus o$, and
for all $h\in(0,1]$,
$|\sigma_{2,h}(P_h)(x,\xi)^{-1}|\leq
C |\xi|_x^{-2}$,\
with $|\xi|^2_x=g_x(\xi,\xi)$ the length of $\xi\in T^*_x X$ with
respect to $g$ and $|.|$ is the operator norm of the matrix of
$\sigma_{2,h}(P_h)(x,\xi)^{-1}$ in any (bounded geometry)
trivialization of $E$.

\begin{lemma}\label{lemma:pos-comm}
Suppose $P_h\in\Diff^2_h(X,E)$ is elliptic and
satisfies \eqref{eq:h-comm-16} for some $c>0$. Suppose also
that $\psi\in e^{-\alpha r}L^2(X,E)$ for all $\alpha$. If
\begin{equation}\label{eq:h-eigenfn}
P_h\psi_h=0,\ \psi_h=e^{r/h}\psi,
\end{equation}
then there exists $R>0$ such that $\psi$ vanishes when $r>R$.
\end{lemma}

\begin{rem}
To simplify notation, we drop the bundle $E$ below. Its presence would not
require any changes, except in the notation.
\end{rem}

\begin{proof}
Let $\Psi(X)$ the algebra of
pseudodifferential operators corresponding to the bounded geometry,
with uniform support, see
\cite[Appendix~1, Definition~3.1-3.2]{Shubin:Spectral}, denoted
by $U\Psi(X)$ there.
The elements of $\Psi^0(X)$
are bounded on $L^2(X)$, and if $A\in\Psi^m(X)$ is elliptic, there is
$B\in\Psi^{-m}(X)$ such that $AB-\Id,BA-\Id\in\Psi^{-\infty}(X)$,
so elliptic regularity statements and estimates work as usual.

We also need the corresponding semiclassical space of operators
$\Psi_h(X)$. These can be defined by modifying the definition of $\Psi(X)$
exactly as if $X$ were compact, i.e.\ defining $\Psi_h^m(X)$ near the
diagonal using the semiclassical quantization
of symbols $a$, and globally as the sums of such operators and elements
of $\Psi_h^{-\infty}(X)$. The latter space consists of operators
with smooth Schwartz kernel that
decays rapidly off the diagonal as $h\to 0$. More precisely, for
$R\in \Psi^{-\infty}_h(X)$ we require that its Schwartz kernel $K$
satisfy $K\in\Cinf((0,1]\times
X\times X)$, that there is $C_R>0$ such that $K(x,y)=0$
if $d(x,y)>C_R$, and for all $N$ there is $C_N>0$ such that for all
$\alpha,\beta$ with $|\alpha|\leq N$, $|\beta|\leq N$, and for all $h\in(0,1]$,
\begin{equation*}
|\pa_x^\alpha \pa_y^\beta K(x,y,h)|\leq C_{N}h^{-n}(1+d(x,y)/h)^{-N},
\end{equation*}
in canonical coordinates, with $n=\dim X$.
All standard properties of semiclassical
ps.d.o's remain valid -- indeed here we only require basic elliptic regularity.
The use of ps.d.o.'s can be eliminated,
if desired, by proving the elliptic regularity estimates directly.

Since $P_h$ is an elliptic family, \eqref{eq:h-eigenfn} and elliptic
regularity give
\begin{equation}\label{eq:h-ell-reg}
\|\psi_h\|_{H^2_h(X)}\leq C_1\|\psi_h\|_{L^2(X)},
\end{equation}
$C_1$ independent of $h\in(0,1]$.
Correspondingly, we do not specify below which Sobolev norms we
are taking. In general,
the letter $C,C'$ will be used denote a constant independent of $h\in(0,1]$,
which may vary from line to line.

We first remark that by the Cauchy-Schwarz inequality, and as
$\|R_j^*\psi_h\|\leq C\|\psi_h\|$, $j=1,2,3,4$,
\begin{equation}\begin{split}\label{eq:R1-R2}
&|\langle h R_1\re P_h\psi_h, \psi_h\rangle|\leq C h\|\psi_h\|\|\re P_h\psi_h\|
\leq Ch\|\psi_h\|^2+C h\|\re P_h\psi_h\|^2,\\
&|\langle h R_2\im P_h \psi_h, \psi_h\rangle|\leq Ch\|\psi_h\|\|\im P_h\psi_h\|
\leq Ch\|\psi_h\|^2+Ch\|\im P_h\psi_h\|^2.
\end{split}\end{equation}
Next,
\begin{equation}\begin{split}\label{eq:R3}
&|\langle \psi_h,h^2R_3\psi_h\rangle|\leq Ch^2\|\psi_h\|^2.
\end{split}\end{equation}
Since $R_4$ is supported in $r\leq r_1$, we can take some
$\chi\in\Cinf(\Real)$ identically $1$ on $(-\infty,3r_1/2)$,
supported in $(-\infty,2r_1)$, and deduce that
\begin{equation*}
|\langle\psi_h,h R_4\psi_h\rangle|=|\langle\chi(r)\psi_h,h R_4\chi(r)
\psi_h\rangle|\leq h\|\chi(r)\psi_h\|^2_{H^1_h(X)}.
\end{equation*}
Now, for $r\leq 2r_1$, $|\psi_h|=e^{r/h}|\psi|\leq
e^{2r_1/h}|\psi|$, with a similar estimate for the semiclassical
derivatives, so
\begin{equation}\begin{split}\label{eq:R4}
|\langle\psi_h,h R_4\psi_h\rangle|&\leq h\|\chi(r)\psi_h\|^2_{H^1_h(X)}\\
&\leq Ch e^{4r_1/h}\|\psi\|^2_{H^1_h(X)}\leq C
h e^{4r_1/h}\|\psi\|^2_{H^1(X)}\leq C' he^{4r_1/h}\|\psi\|^2.
\end{split}\end{equation}
Hence, we deduce from \eqref{eq:h-comm-8} (with $P_\alpha$ replaced by $P_h$)
and \eqref{eq:h-comm-16} that
\begin{equation}\begin{split}\label{eq:h-comm-32}
0\geq (1-Ch)\|\re P_h\psi_h\|^2+(1-Ch)\|\im P_h\psi_h\|^2
&+h(c-Ch)\|\psi_h\|^2\\
&-Che^{4r_1/h}\|\psi\|^2.
\end{split}\end{equation}
Dropping the first two (positive) terms on the right hand side, we conclude
that there exists $h_0>0$ such that for $h\in(0,h_0)$,
\begin{equation}\label{eq:h-comm-64}
Ch e^{4r_1/h}\|\psi\|^2
\geq h\frac{c}{2}\|\psi_h\|^2.
\end{equation}

Now suppose that $R>2r_1$ and
$\supp\psi\cap\{
r\geq R\}$ is non-empty.
Since $e^{2r/h}\geq e^{2R/h}$ for $r\geq R$, we deduce
that
\begin{equation*}
\|\psi_h\|^2\geq  C' e^{2R/h},\ C'
=\|\psi\|^2_{r\geq R}>0.
\end{equation*}
Thus, we conclude from \eqref{eq:h-comm-64} that
\begin{equation}\label{eq:h-lim}
C\|\psi\|^2\geq \frac{c}{2}\,C'
e^{2(R-2r_1)/h}.
\end{equation}
But letting $h\to 0$,
the right hand side goes to $+\infty$, providing a contradiction.

Thus, $\psi$ vanishes for $r\geq R$.
\end{proof}

The proof of Theorem~\ref{thm:Lap} is finished since if $\psi$ vanishes on an
open set, it vanishes everywhere on $X$
by the usual Carleman-type unique continuation theorem
\cite[Theorem~17.2.1]{Hor}.

In fact, it is straightforward to strengthen Lemma~\ref{lemma:pos-comm}
and allow
$P\psi$ to be compactly supported. The following lemma thus completes
the proof of Theorem~\ref{thm:perturb}:

\begin{lemma}
Suppose $P\in\Diff^2(X;E)$ is elliptic, $\psi\in e^{-\alpha r}L^2(X,E)$
for all $\alpha$, and there is $r_0>0$ such that $P\psi=0$ for $r>r_0$.
Let $P_h=e^{r/h}h^2P e^{-r/h}$, and suppose that $P_h$
satisfies \eqref{eq:h-comm-16} for some $c>0$.
Then there exists $R>0$ such that $\psi$ vanishes when $r>R$.
\end{lemma}

\begin{rem}\label{rem:loc}
Note that in this formulation, if $X$ is replaced by a manifold with several
ends, one of which is of the product form eluded to in the
introduction, our theorem holds locally on this end. That is,
if $P\psi$ vanishes on this end and $\psi$ has superexponential decay there,
then $\psi$ vanishes on the end --- hence globally by the standard
unique continuation theorem if $P\psi$ is identically zero.
To prove this, we merely multiply by a cutoff function supported on this
end, and apply the lemma to the resulting inhomogeneous problem.
\end{rem}

\begin{proof}
The elliptic regularity estimate now becomes
\begin{equation}\label{eq:h-ell-reg-p}
\|\psi_h\|_{H^2_h(X)}\leq C_1(\|\psi_h\|_{L^2(X)}+\|P_h\psi_h\|_{L^2(X)}),
\end{equation}
$C_1$ independent of $h\in(0,1]$, and we need to keep track of the
second term on the right hand side.

Correspondingly,
$\|R_j^*\psi_h\|\leq C(\|\psi_h\|+\|P_h\psi_h\|)$, $j=1,2,3,4$.
Thus, on the right hand side of \eqref{eq:R1-R2}, we need to add
$Ch \|P_h\psi_h\|^2$, resp.\ $Ch \|P_h\psi_h\|^2$, while
on the right hand side of \eqref{eq:R3} we need to add $Ch^2\|P_h\psi_h\|^2$.
Similarly, we need to add $C'he^{4r_1/h}\|P\psi\|^2$ to the right hand side
of \eqref{eq:R4}.
Thus, \eqref{eq:h-comm-32} becomes
\begin{equation*}\begin{split}
(1+Ch)
\|P_h\psi_h\|^2
\geq (1-Ch)&\|\re P_h\psi_h\|^2+(1-Ch)\|\im P_h\psi_h\|^2\\
&+h(c-Ch)\|\psi_h\|^2
-Che^{4r_1/h}(\|\psi\|^2+\|P\psi\|^2).
\end{split}\end{equation*}
Since $P_h\psi_h=e^{r/h}h^2P\psi$, we have $\|P_h\psi_h\|\leq e^{r_0/h}h^2
\|P\psi\|$. Let $r_2=\max(r_0,r_1)$.
Thus, there exists $h_0>0$ such that for $h\in(0,h_0)$,
\begin{equation}\label{eq:h-comm-64-p}
2e^{4r_2/h}\|P\psi\|^2+Ch e^{4r_1/h}\|\psi\|^2
\geq h\frac{c}{2}\|\psi_h\|^2.
\end{equation}
Taking $R>2r_2$,
the proof is now finished as in Lemma~\ref{lemma:pos-comm}, for
\eqref{eq:h-lim} becomes
\begin{equation*}
2\|P\psi\|^2+Ch\|\psi\|^2\geq \frac{c}{2}\,C'
he^{2(R-2r_2)/h},
\end{equation*}
and the right hand side still goes to $+\infty$, while the left hand side
is bounded as $h\to 0$.
\end{proof}

\def\cprime{$'$} \def\cprime{$'$}

\end{document}